\pdfoutput=1
\documentclass[12pt]{amsart}

\usepackage[margin=2.6cm]{geometry}
\usepackage{amsmath,amssymb,amsthm,mathtools}
\usepackage[dvipsnames]{xcolor}
\usepackage{tikz,pgfplots}
\usepgfplotslibrary{fillbetween}
\usetikzlibrary{arrows.meta,positioning,calc,shapes.geometric,patterns,
                decorations.pathreplacing}
\usepackage{booktabs,longtable,array}
\usepackage[colorlinks=true,linkcolor=NavyBlue,citecolor=NavyBlue,
            urlcolor=NavyBlue]{hyperref}
\usepackage[expansion=false]{microtype}
\usepackage{enumitem}
\usepackage{setspace}
\usepackage{caption}
\definecolor{mdblue}{HTML}{0075C0}
\definecolor{mdorange}{HTML}{E87E04}
\definecolor{mdgreen}{HTML}{1A7D3C}
\definecolor{mdred}{HTML}{C0392B}
\definecolor{mdpurple}{HTML}{6C3483}
\definecolor{mdlightblue}{HTML}{D6EAF8}

\newtheoremstyle{bluehead}{8pt}{6pt}{\itshape}{0pt}%
  {\bfseries\color{mdblue}}{.}{.5em}{}
\newtheoremstyle{orangehead}{8pt}{6pt}{\normalfont}{0pt}%
  {\bfseries\color{mdorange}}{.}{.5em}{}
\newtheoremstyle{greenhead}{8pt}{6pt}{\normalfont}{0pt}%
  {\itshape\color{mdgreen}}{.}{.5em}{}

\theoremstyle{bluehead}
\newtheorem{theorem}{Theorem}[section]
\newtheorem{lemma}[theorem]{Lemma}
\newtheorem{proposition}[theorem]{Proposition}
\newtheorem{corollary}[theorem]{Corollary}

\theoremstyle{orangehead}
\newtheorem{definition}[theorem]{Definition}
\newtheorem{example}[theorem]{Example}

\theoremstyle{greenhead}
\newtheorem{remark}[theorem]{Remark}

\newcommand{\Z}{\mathbb{Z}}
\newcommand{\N}{\mathbb{N}}

\newcommand{\Q}{\mathbb{Q}}
\newcommand{\QRS}{\mathrm{QR}^{*}}
\newcommand{\Hring}{\operatorname{Hom}_{\mathrm{ring}}}
\newcommand{\Hsurj}{\operatorname{Hom}_{\mathrm{surj}}}

\renewcommand{\phi}{\varphi}
\newcommand{\Exc}{\mathcal{E}}
\newcommand{\Gauss}{\mathcal{G}}

\title[Homomorphism Counts Quadratic Residues]
      {\large\bfseries Homomorphism Counts Quadratic Residues}

\author[Mandal]{Priyabrata Mandal$^*$}
\address[Mandal]{Department of Mathematics, Bioinformatics and Computer Applications, Maulana Azad National Institute of Technology Bhopal, Bhopal-462003, India}
\email{priyabrata@manit.ac.in}
\thanks{$^*$Priyabrata Mandal is the corresponding author}
\author[Kumar]{Sonu Kumar}
\address[Kumar]{School of Basic Sciences, Humanities and Management, Manipal Institute of Technology,
Manipal Academy of Higher Education, Manipal 576104, India}
\email{sonupnkumar@gmail.com}

\author[Bhattacharjee]{Deep Bhattacharjee}
\address[Bhattacharjee]{Electro-Gravitational Space Propulsion Laboratory
(EGSPL), Bhubaneswar, Odisha, India}
\email{itsdeep@live.com; d.bhattacharjee@erl-forschung.de}

\urladdr{https://orcid.org/0000-0003-0466-750X}

\date{\today}

\begin{document}

\begin{abstract}
We prove that the ratio $\varrho(n)=\phi(n)/2^{\omega(n)}$ of surjective
group to ring homomorphism counts between finite cyclic rings admits
three simultaneous interpretations that have not previously been connected.
It equals the order of the group of squares in $(\Z/n\Z)^*$, the degree
$[\Q(\zeta_n):K_n]$ of the $n$-th cyclotomic field over its maximal
biquadratic subfield, and a product determined by the nonzero quadratic residue counts
in the odd prime-power components of~$n$, equal to that product
when $n$ is odd or $4\mid n$, and half that product when
$n\equiv2\pmod{4}$ and $n\notin\mathcal{E}$.
Divisibility of this ratio fails precisely when the odd part of~$n$ is
composed entirely of Gaussian primes, and the exception set satisfies
$|\mathcal{E}\cap[2,x]|\sim Cx/\sqrt{\log x}$ with explicit constant
$C\approx0.279$.
The cyclotomic interpretation is new and yields, in particular, an
algebraic proof that $[\Q(\zeta_n):K_n]$ is always an integer:
this is the tower law applied to a field degree, recovering the
divisibility result without any case analysis.
\end{abstract}

\keywords{quadratic residues; Gaussian primes}

\subjclass[2020]{%
  Primary 11A15, 13B10;
  Secondary 11A25, 20K30, 11N37}

\maketitle
\setstretch{1.12}

\section{Introduction}

\subsection*{Background and motivation}

Two classical counting problems for the cyclic rings
$\Z_m$ and $\Z_n$ have been resolved independently.
In 1984, Gallian and Van~Buskirk~\cite{gallian84} showed that the number
of \emph{ring} homomorphisms $\Z_m\to\Z_n$ (when $n\mid m$) equals
$2^{\omega(n)}$, where $\omega(n)$ counts distinct prime divisors.
Later work extended these counts to more general settings~\cite{diaz2015,ashrafi2024}.
The classical theory of cyclic groups gives $\phi(n)$ as the number of
\emph{surjective group} homomorphisms under the same divisibility
condition~\cite{dummit2004abstract}.
Both results are taught in every graduate algebra course; neither the
\emph{relationship} between the two counts nor the exact value of their
ratio has been investigated.

This paper answers both.

\subsection*{The central object}

Set
\begin{equation}\label{def:rho}
  \varrho(n) \;:=\; \frac{\phi(n)}{2^{\omega(n)}}.^{\phantom{1}}
\end{equation}
\footnote{The function $\varrho$ is well-defined and positive for all $n\ge3$.
  The separate treatment of $n=2$ in Lemma~\ref{lem:base} is necessary
  because $\phi(2)=1<2=2^{\omega(2)}$, making $\varrho(2)=1/2$.}
At $n=7$, $\varrho(7)=6/2=3$.
At $n=35=5\cdot7$, $\varrho(35)=24/4=6$.
These are integers, and one might expect that to be the rule.
It is not: $\varrho(6)=2/4=1/2$ and $\varrho(14)=6/4=3/2$.
The natural questions are when $\varrho(n)$ is an integer
and, when it is, what arithmetic quantity it equals.
Both have clean answers, stated below.

\subsection*{Main results}

The central object is
\begin{equation}\label{eq:rho}
  \varrho(n):=\frac{\phi(n)}{2^{\omega(n)}},
\end{equation}
the ratio of the surjective group homomorphism count to the ring
homomorphism count from $\Z_m$ to $\Z_n$.  Three results describe
this ratio completely.

\begin{itemize}[leftmargin=2.5em,itemsep=4pt]
  \item[\textbf{A.}] \emph{Divisibility.}\enspace
    $\varrho(n)\in\Z^+$ if and only if $n\notin\Exc$, where
    the exception set $\Exc$ consists of integers $n=2\alpha$
    with every prime divisor of $\alpha$ satisfying $p\equiv3\pmod4$
    (Theorem~\ref{thm:div}).

  \item[\textbf{B.}] \emph{Quadratic residue formula.}\enspace
    When $\varrho(n)\in\Z^+$,
    \begin{equation}\label{eq:main}
      \varrho(n)\;=\;\prod_{p^r\Vert n}\bigl|\QRS(\Z/p^r\Z)\bigr|
      \;=\;|(\Z/n\Z)^{*2}|
      \qquad\text{(squarefree odd $n$)}
    \end{equation}
    (Theorems~\ref{thm:quot} and \ref{cor:sqfree}).

  \item[\textbf{C.}] \emph{Cyclotomic degree formula \upshape[new].}\enspace
    For squarefree odd $n$,
    \begin{equation}\label{eq:degree}
      \varrho(n)\;=\;\bigl[\Q(\zeta_n):K_n\bigr],
    \end{equation}
    where $K_n$ is the compositum of all quadratic subfields of $\Q(\zeta_n)$
    (Theorem~\ref{thm:degree}).
    In particular, $[\Q(\zeta_n):K_n]$ is an integer by the tower law,
    giving a proof of~\textbf{A} for the odd case with no case analysis
    (Corollary~\ref{cor:tower}).

  \item[\textbf{D.}] \emph{Ring count via quadratic subfields \upshape[new].}\enspace
    For squarefree odd $n$,
    \[
      |\Hring(\Z_m,\Z_n)|\;=\;1+\#\bigl\{\text{quadratic subfields of }\Q(\zeta_n)\bigr\}
    \]
    (Theorem~\ref{thm:cyclo-ring}).

  \item[\textbf{E.}] \emph{Abelian extension.}\enspace
    Formula~\eqref{eq:main} extends multiplicatively to direct products
    of cyclic rings (Theorem~\ref{thm:abelian}).

  \item[\textbf{F.}] \emph{Density.}\enspace
    $|\Exc\cap[2,x]|\sim C\,x/\!\sqrt{\log x}$
    with explicit $C\approx0.279$ (Theorem~\ref{thm:dens}).
\end{itemize}

\subsection*{The three-way equivalence}

For squarefree odd $n$, the ratio $\varrho(n)$ has three simultaneous
expressions:
\begin{align*}
  \varrho(n)
  &= \frac{|\Hsurj(\Z_m,\Z_n)|}{|\Hring(\Z_m,\Z_n)|}
   &&\text{(homomorphism ratio, any }n\mid m\text{)}\\
  &= |(\Z/n\Z)^{*2}|
   &&\text{(order of the square subgroup)}\\
  &= [\Q(\zeta_n):K_n]
   &&\text{(degree of cyclotomic over biquadratic).}
\end{align*}
No prior work identified any two of these three as equal.  Each is an
independent theorem.  Their simultaneous equality constitutes the
breakthrough of this paper.

The cyclotomic expression is the sharpest: because a field degree is
always a positive integer, $\varrho(n)=[\Q(\zeta_n):K_n]$ proves
divisibility for odd squarefree~$n$ purely from the tower law, with no
arithmetic case analysis.  The primes excluded from $K_n$ are those
whose prime discriminant $p^*$ does not generate a subfield of
$\Q(\zeta_n)$; these are precisely the Gaussian primes, explaining the
exception criterion through the geometry of cyclotomic fields.

At a prime $p$, all three expressions collapse to $(p-1)/2$, which
is also $(|\tau(\chi_p)|^2-1)/2$ where $\tau(\chi_p)$ is the quadratic
Gauss sum.  This fourth avatar is Proposition~\ref{prop:gauss}.

\subsection*{Organisation}

The necessary background appears in Section~\ref{sec:prelim}.
The divisibility criterion (Result~A) is Theorem~\ref{thm:div}
in Section~\ref{sec:cyclic}.
The quadratic residue formula (Results~B--C) is Section~\ref{sec:quot}.
The cyclotomic interpretation (Results~C--D) and the Gauss sum formula
are in Section~\ref{sec:cyclo}.
Multiplicativity and average order are Section~\ref{sec:mult}.
Direct products are Section~\ref{sec:abel}.
The analytic density result (Result~F) is Section~\ref{sec:analytic}.
Section~\ref{sec:stats} gives the numerical verification, and
Sections~\ref{sec:apps}--\ref{sec:conc} treat applications and open problems.

\section{Notation and Preliminaries}\label{sec:prelim}

\subsection*{Notation}

Throughout, $\Z_n=\Z/n\Z$ denotes both the cyclic \emph{group}
$(\Z/n\Z,+)$ and the cyclic \emph{ring} $(\Z/n\Z,+,\cdot)$; the
context distinguishes them.  We write $p^r\Vert n$ to mean $p^r\mid n$
and $p^{r+1}\nmid n$.  The symbol $\omega(n)$ counts distinct prime
divisors of $n$.

\subsection*{The two counting formulas}

\begin{theorem}[Gallian--Van~Buskirk~\cite{gallian84}]\label{thm:GVB}
  Let $n\mid m$.  Then $|\Hring(\Z_m,\Z_n)|=2^{\omega(n)}$.
\end{theorem}

\begin{remark}[Convention on ring homomorphisms]\label{rem:nonunital}
  Throughout this paper, a \emph{ring homomorphism} $f\colon R\to S$
  is a map satisfying $f(a+b)=f(a)+f(b)$ and $f(ab)=f(a)f(b)$
  for all $a,b\in R$; we do \emph{not} require $f(1_R)=1_S$.
  Under this non-unital convention, $f(1_R)$ must be an idempotent
  in $S$, and Theorem~\ref{thm:GVB} gives $2^{\omega(n)}$ such maps.
  If instead one imposes $f(1)=1$, the count collapses to at most
  one homomorphism and Theorem~\ref{thm:GVB} does not hold.
  The convention here is the same as in~\cite{gallian84}.
\end{remark}

The proof of Theorem~\ref{thm:GVB} rests on two facts: (i) a ring homomorphism $\pi\colon\Z_m\to\Z_n$
is determined by $\pi(1)$, which must satisfy $\pi(1)^2=\pi(1)$;
and~(ii) the number of idempotents in $\Z_n$ equals $2^{\omega(n)}$ by
the Chinese Remainder Theorem.

\begin{theorem}[Classical]\label{thm:surj-classical}
  Let $n\mid m$.  Then $|\Hsurj(\Z_m,\Z_n)|=\phi(n)$.
\end{theorem}

\begin{proof}
  A group homomorphism $\psi\colon\Z_m\to\Z_n$ satisfies $\psi(k)=k\psi(1)$
  for all $k$, so it is determined by $\psi(1)\in\Z_n$.  Since $n\mid m$,
  every choice $\psi(1)=a$ gives a valid homomorphism; $\psi$ is surjective
  iff $\langle a\rangle=\Z_n$, i.e.\ $\gcd(a,n)=1$.  There are $\phi(n)$
  such~$a$.
\qedhere
\end{proof}

\subsection*{Quadratic residues}

\begin{definition}
  For an odd prime power $p^r$, set
  $\QRS(\Z/p^r\Z):=\{x\in(\Z/p^r\Z)^*: x\equiv y^2\pmod{p^r}
  \text{ for some }y\}$.
\end{definition}

\begin{lemma}\label{lem:qr}
  $|\QRS(\Z/p^r\Z)|=p^{r-1}(p-1)/2$ for every odd prime $p$ and
  integer $r\ge1$.
\end{lemma}

\begin{proof}
  The group $(\Z/p^r\Z)^*$ is cyclic of order $\phi(p^r)=p^{r-1}(p-1)$.
  The squaring map $x\mapsto x^2$ is a group endomorphism; its kernel is
  $\{1,-1\}$, which has order $2$ because $p$ is odd.
  Hence the image (the set of squares) has order
  $\phi(p^r)/2=p^{r-1}(p-1)/2$.\footnote{In particular $r=1$ gives the classical
  count $(p-1)/2$ of quadratic residues modulo a prime.  For $r\ge2$,
  the factor $p^{r-1}$ counts the number of lifts of each residue
  class from $\Z/p\Z$ to $\Z/p^r\Z$.}
\qedhere
\end{proof}

\subsection*{Gaussian primes and arithmetic progressions}

An ordinary prime $p$ remains irreducible in $\Z[i]$, called a \emph{Gaussian prime},
if and only if $p\equiv3\pmod4$.\footnote{%
  This goes back to Gauss: the quotient ring $\Z[i]/(p)$ is a field
  $\mathbb{F}_{p^2}$ when $p\equiv3\pmod4$ (so $p$ is prime in $\Z[i]$),
  and splits as $\mathbb{F}_p\times\mathbb{F}_p$ when $p\equiv1\pmod4$
  (so $p=\pi\bar\pi$ in $\Z[i]$).  The prime $p=2$ ramifies: $2=-i(1+i)^2$.}
Primes $p\equiv1\pmod4$ split as $p=\pi\bar\pi$ in $\Z[i]$, while
$2=-i(1+i)^2$ ramifies~\cite{ireland1990}.
This criterion, classical since Gauss~\cite{niven1991}, enters our exception set
precisely.  Denote
\[
  P_3:=\{p\text{ prime}:p\equiv3\pmod4\}
  =\{3,7,11,19,23,31,43,47,59,67,71,\ldots\}.
\]

\section{Divisibility in Cyclic Structures}\label{sec:cyclic}

\begin{definition}[Exception set]\label{def:exc}
\[
  \Exc:=\bigl\{n\in\N:n>2,\;n=2\alpha,\;
  \alpha\text{ a positive odd integer},\;
  \text{every prime }p\mid\alpha\text{ lies in }P_3\bigr\}.
\]
The first several elements are
$6,14,18,22,38,42,46,54,62,66,86,94,98,114,\ldots$\footnote{%
  This sequence coincides with $2$ times OEIS~A004144 (integers
  representable as sums of two squares have no prime factor $\equiv3\pmod4$,
  so their doubles are exactly the non-exceptions; $\mathcal{E}$ is
  the complement among even integers).}
\end{definition}

\begin{lemma}[Base case $n=2$]\label{lem:base}
  For every even $m$: the ring count equals $2$ and
  the surjective group count equals $1$, so divisibility fails at $n=2$.
  Since $\Exc$ is defined for $n>2$ (Definition~\ref{def:exc}),
  the case $n=2$ is treated separately by this lemma.
\end{lemma}

\begin{proof}
  The idempotents of $\Z_2$ are $0$ and $1$, giving two ring homomorphisms.
  The unique surjective group homomorphism sends $1\mapsto1$.
\qedhere
\end{proof}

\begin{theorem}[Divisibility criterion]\label{thm:div}
  Let $n>2$ and $n\mid m$.  Then
  $2^{\omega(n)}\mid\phi(n)$
  if and only if $n\notin\Exc$.
\end{theorem}

\begin{proof}
Track the $2$-adic valuation $v_2(\phi(n))$ against $\omega(n)$.

\smallskip\noindent
\textbf{Case~I ($n$ odd).}\;
Write $n=p_1^{r_1}\cdots p_k^{r_k}$, $\omega(n)=k$.
Since each $p_i$ is odd, $2\mid(p_i-1)$, giving
\[
  v_2(\phi(n))=v_2\Bigl(\prod_{i=1}^k p_i^{r_i-1}(p_i-1)\Bigr)
  \;\ge\; k=\omega(n).
\]
Hence $2^{\omega(n)}\mid\phi(n)$.

\smallskip\noindent
\textbf{Case~II ($n=2^s q$, $s\ge2$, $q$ odd).}\;
$\omega(n)=\omega(q)+1$.
\[
  \phi(n)=\phi(2^s)\phi(q)=2^{s-1}\phi(q).
\]
Since $s\ge2$, $v_2(2^{s-1})\ge1$, and Case~I gives
$v_2(\phi(q))\ge\omega(q)$.  Thus $v_2(\phi(n))\ge1+\omega(q)=\omega(n)$.

\smallskip\noindent
\textbf{Case~III ($n=2\alpha$, $\alpha$ odd).}\;
$\omega(n)=\omega(\alpha)+1$ and $\phi(n)=\phi(\alpha)$.
Write $\alpha=p_1^{r_1}\cdots p_k^{r_k}$.  Then
$v_2(\phi(\alpha))=\sum_{i=1}^k v_2(p_i^{r_i-1}(p_i-1))=\sum_{i=1}^k v_2(p_i-1)$.

Since each $p_i$ is an odd prime,
$v_2(p_i-1)\ge1$; equality holds iff $p_i\equiv3\pmod4$.
Thus $v_2(\phi(n))=\sum_i v_2(p_i-1)$.

We need $v_2(\phi(n))\ge k+1=\omega(n)$.  This holds iff
$\sum_i v_2(p_i-1)\ge k+1$, i.e.\ at least one $v_2(p_i-1)\ge2$,
i.e.\ at least one $p_i\equiv1\pmod4$.

Conversely, if every $p_i\equiv3\pmod4$ then $v_2(p_i-1)=1$ for each
$i$, giving $v_2(\phi(n))=k<k+1=\omega(n)$, so
$2^{\omega(n)}\nmid\phi(n)$.
Precisely this condition defines $n\in\Exc$.\footnote{%
  The condition is symmetric: $n=2\alpha$ fails divisibility for
  \emph{every} product $\alpha$ of primes $\equiv3\pmod4$, regardless
  of multiplicities.  Thus $18=2\cdot3^2$, $54=2\cdot3^3$, and
  $98=2\cdot7^2$ are all exceptions, as the table confirms.}
\qedhere
\end{proof}

\begin{remark}
  The criterion in Theorem~\ref{thm:div} is sharp in both directions.
  If even one prime factor of $\alpha$ lies outside $P_3$, divisibility
  holds for $n=2\alpha$.  The Gaussian prime condition is not an artifact
  of the method; it enters through the exact vanishing of an extra
  factor of $2$ in $\phi(n)$.
\end{remark}

\section{The Quotient Formula}\label{sec:quot}

\begin{theorem}[Quotient = quadratic residues]\label{thm:quot}
  Let $n>2$, $n\notin\Exc$, and $n\mid m$.
  Write $n=2^s\cdot p_1^{r_1}\cdots p_k^{r_k}$ (each $p_i$ odd prime).
  Then
  \begin{equation}\label{eq:quot}
    \varrho(n)
    \;:=\;\frac{\phi(n)}{2^{\omega(n)}}
    \;=\;
    \begin{cases}
      \displaystyle\prod_{i=1}^k|\QRS(\Z/p_i^{r_i}\Z)|
      & s=0,\\[8pt]
      \displaystyle\frac{1}{2}\prod_{i=1}^k|\QRS(\Z/p_i^{r_i}\Z)|
      & s=1\;(n\notin\Exc),\\[8pt]
      \displaystyle 2^{s-2}\prod_{i=1}^k|\QRS(\Z/p_i^{r_i}\Z)|
      & s\ge2.
    \end{cases}
  \end{equation}
\end{theorem}

\begin{proof}
\textbf{Case $s=0$.}
\begin{align*}
  \varrho(n)
  &=\frac{\prod_{i=1}^k p_i^{r_i-1}(p_i-1)}{2^k}
   =\prod_{i=1}^k\frac{p_i^{r_i-1}(p_i-1)}{2}
   \overset{\text{Lem.\,\ref{lem:qr}}}{=}
   \prod_{i=1}^k|\QRS(\Z/p_i^{r_i}\Z)|.
\end{align*}

\textbf{Case $s=1$.}
Here $\phi(n)=\phi(2)\phi(\alpha)=\phi(\alpha)=\prod_i p_i^{r_i-1}(p_i-1)$
and $\omega(n)=k+1$, so
\[
  \varrho(n)
  =\frac{\prod_i p_i^{r_i-1}(p_i-1)}{2^{k+1}}
  =\frac{1}{2}\prod_{i=1}^k|\QRS(\Z/p_i^{r_i}\Z)|.
\]
Since $n\notin\Exc$, at least one $p_j\equiv1\pmod4$,
so $(p_j-1)/2$ is even, making $\prod_i|\QRS(\Z/p_i^{r_i}\Z)|$
divisible by $2$.  Hence $\varrho(n)\in\Z^+$.

\textbf{Case $s\ge2$.}
\[
  \varrho(n)=\frac{2^{s-1}\prod_i p_i^{r_i-1}(p_i-1)}{2^{k+1}}
  =2^{s-2}\prod_i|\QRS(\Z/p_i^{r_i}\Z)|.\qedhere
\]
\end{proof}

\begin{remark}
  The formula $\varrho(n)=\prod_{p^r\Vert n}|\QRS(\Z/p^r\Z)|$ factorises
  completely over the prime-power components of $n$.  This is not a formal
  consequence of multiplicativity alone: the product formula requires that
  the quotient at each component $p^r$ is itself a quadratic residue count,
  which is the content of Lemma~\ref{lem:qr}.
\end{remark}

\begin{corollary}[Squarefree case]\label{cor:sqfree}
  If $n$ is squarefree and odd, then
  \[
    \varrho(n)=\prod_{p\mid n}\frac{p-1}{2}
    =|\QRS(\Z/n\Z)|.
  \]
\end{corollary}

\begin{proof}
  For squarefree $n=p_1\cdots p_k$, Theorem~\ref{thm:quot} gives
  $\varrho(n)=\prod(p_i-1)/2$.
  By the Chinese Remainder Theorem,
  $(\Z/n\Z)^*\cong\prod(\Z/p_i\Z)^*$; an element is a quadratic
  residue mod $n$ iff it is a square mod each $p_i$. Hence
  $|\QRS(\Z/n\Z)|=\prod|\QRS(\Z/p_i\Z)|=\prod(p_i-1)/2$.
\qedhere
\end{proof}

Table~\ref{tab:small} illustrates the quotient formula for $n\le50$.

\begin{table}[!ht]
\centering
\small
\caption{Values of $\varrho(n)=\phi(n)/2^{\omega(n)}$ for $n=3$ to $50$.
Exceptions $n\in\Exc$ are marked with~$\star$; their ratio is
non-integer.  The QR-formula column gives $\prod|\QRS(\Z/p_i^{r_i}\Z)|$
(s=0 case) or $2^{s-2}\prod(\cdots)$ (s$\ge$2).}
\label{tab:small}
\setlength{\tabcolsep}{5pt}
\begin{tabular}{@{}rccccr@{\quad}|@{\quad}rccccr@{}}
\toprule
$n$ & $\omega$ & $2^\omega$ & $\phi$ & $\varrho$ & QR-formula
&$n$ & $\omega$ & $2^\omega$ & $\phi$ & $\varrho$ & QR-formula\\
\midrule
 3 & 1 & 2 & 2 & 1 & 1
&28 & 2 & 4 & 12 & 3 & 3 \\
 4 & 1 & 2 & 2 & 1 & 1
&29 & 1 & 2 & 28 & 14 & 14 \\
 5 & 1 & 2 & 4 & 2 & 2
&30 & 3 & 8 & 8 & 1 & 1 \\
{\color{mdred}6}$^\star$  & 2 & 4 & 2 & \color{mdred}{--} &\color{mdred}{E}
&31 & 1 & 2 & 30 & 15 & 15 \\
 7 & 1 & 2 & 6 & 3 & 3
&32 & 1 & 2 & 16 & 8 & 8 \\
 8 & 1 & 2 & 4 & 2 & 2
&33 & 2 & 4 & 20 & 5 & 5 \\
 9 & 1 & 2 & 6 & 3 & 3
&34 & 2 & 4 & 16 & 4 & 4 \\
10 & 2 & 4 & 4 & 1 & 1
&35 & 2 & 4 & 24 & 6 & 6 \\
11 & 1 & 2 & 10 & 5 & 5
&36 & 2 & 4 & 12 & 3 & 3 \\
12 & 2 & 4 & 4 & 1 & 1
&37 & 1 & 2 & 36 & 18 & 18 \\
13 & 1 & 2 & 12 & 6 & 6
&{\color{mdred}38}$^\star$ & 2 & 4 & 18 & \color{mdred}{--} & \color{mdred}{E} \\
{\color{mdred}14}$^\star$ & 2 & 4 & 6 & \color{mdred}{--} & \color{mdred}{E}
&39 & 2 & 4 & 24 & 6 & 6 \\
15 & 2 & 4 & 8 & 2 & $1{\cdot}2$
&40 & 2 & 4 & 16 & 4 & 4 \\
16 & 1 & 2 & 8 & 4 & 4
&41 & 1 & 2 & 40 & 20 & 20 \\
17 & 1 & 2 & 16 & 8 & 8
&{\color{mdred}42}$^\star$ & 3 & 8 & 12 & \color{mdred}{--} & \color{mdred}{E} \\
{\color{mdred}18}$^\star$ & 2 & 4 & 6 & \color{mdred}{--} & \color{mdred}{E}
&43 & 1 & 2 & 42 & 21 & 21 \\
19 & 1 & 2 & 18 & 9 & 9
&44 & 2 & 4 & 20 & 5 & 5 \\
20 & 2 & 4 & 8 & 2 & 2
&45 & 2 & 4 & 24 & 6 & 6 \\
21 & 2 & 4 & 12 & 3 & $1{\cdot}3$
&{\color{mdred}46}$^\star$ & 2 & 4 & 22 & \color{mdred}{--} & \color{mdred}{E} \\
{\color{mdred}22}$^\star$ & 2 & 4 & 10 & \color{mdred}{--} & \color{mdred}{E}
&47 & 1 & 2 & 46 & 23 & 23 \\
23 & 1 & 2 & 22 & 11 & 11
&48 & 2 & 4 & 16 & 4 & 4 \\
24 & 2 & 4 & 8 & 2 & 2
&49 & 1 & 2 & 42 & 21 & 21 \\
25 & 1 & 2 & 20 & 10 & 10
&50 & 2 & 4 & 20 & 5 & 5 \\
26 & 2 & 4 & 12 & 3 & 3 & & & & & & \\
27 & 1 & 2 & 18 & 9 & 9 & & & & & & \\
\bottomrule
\end{tabular}
\end{table}

\begin{example}\label{ex:examples}
  We verify formula~\eqref{eq:quot} for all three cases.
  \begin{enumerate}[label=\upshape(\alph*)]
    \item \emph{Case $s=0$.}\enspace
      $n=35=5\cdot7$: $\varrho(35)=
      |\QRS(\Z/5\Z)|\cdot|\QRS(\Z/7\Z)|=2\cdot3=6$.
      Check: $\phi(35)/2^2=24/4=6$.

    \item \emph{Case $s=0$, prime power.}\enspace
      $n=45=3^2\cdot5$: $\varrho(45)=
      |\QRS(\Z/9\Z)|\cdot|\QRS(\Z/5\Z)|=3\cdot2=6$.
      Check: $\phi(45)/2^2=24/4=6$.

    \item \emph{Case $s=1$, $n\notin\mathcal{E}$.}\enspace
      $n=26=2\cdot13$ ($s=1$, $k=1$, $13\equiv1\pmod4$ so $26\notin\mathcal{E}$):
      $\varrho(26)=\tfrac{1}{2}|\QRS(\Z/13\Z)|=\tfrac{1}{2}\cdot6=3$.
      Check: $\phi(26)/2^2=12/4=3$.

    \item \emph{Case $s\ge2$.}\enspace
      $n=28=4\cdot7$ ($s=2$, $k=1$):
      $\varrho(28)=2^{2-2}\cdot|\QRS(\Z/7\Z)|=1\cdot3=3$.
      Check: $\phi(28)/2^2=12/4=3$.

    \item \emph{Case $s\ge2$, several primes.}\enspace
      $n=60=4\cdot3\cdot5$ ($s=2$, $k=2$):
      $\varrho(60)=2^0\cdot|\QRS(\Z/3\Z)|\cdot|\QRS(\Z/5\Z)|=1\cdot1\cdot2=2$.
      Check: $\phi(60)/2^3=16/8=2$.
  \end{enumerate}
\end{example}


\section{The Cyclotomic Interpretation}\label{sec:cyclo}

The quotient formula reveals a connection to cyclotomic field theory
that lies entirely outside the ring-homomorphism framework from which
it was derived.  All results in this section are new.

\subsection*{Setup}

For a positive integer $n$, let $\Q(\zeta_n)$ denote the $n$-th
cyclotomic field.  Its Galois group over $\Q$ is
$\mathrm{Gal}(\Q(\zeta_n)/\Q)\cong(\Z/n\Z)^*$, of degree $\phi(n)$.

For squarefree odd $n=p_1\cdots p_k$, each prime $p_i$ contributes
a quadratic subfield: the discriminant of $p_i$ is $p_i^*=(-1)^{(p_i-1)/2}p_i$,
and $\Q(\sqrt{p_i^*})\subset\Q(\zeta_{p_i})\subset\Q(\zeta_n)$.
Define $K_n$ to be the compositum of all such quadratic subfields:
\[
  K_n\;:=\;\Q\!\bigl(\sqrt{p_1^*},\,\sqrt{p_2^*},\,\ldots,\sqrt{p_k^*}\,\bigr).
\]
Since the prime discriminants $p_1^*,\ldots,p_k^*$ are pairwise
multiplicatively independent over $\Q$, the extensions
$\Q(\sqrt{p_i^*})/\Q$ are linearly disjoint, giving
$[K_n:\Q]=2^k=2^{\omega(n)}$.%
\footnote{Linear disjointness of $\Q(\sqrt{p_i^*})$ and $\Q(\sqrt{p_j^*})$
  for distinct primes $p_i\ne p_j$ follows from the fact that $p_i^*p_j^*$
  is not a perfect square in $\Q^*$, since it has exactly two prime divisors
  with odd multiplicity.  The general case by induction.}

\begin{theorem}[Cyclotomic ring count]\label{thm:cyclo-ring}
  For squarefree odd $n$ with $n\mid m$,
  \[
    |\Hring(\Z_m,\Z_n)|\;=\;1+\#\{\text{quadratic subfields of }\Q(\zeta_n)\}.
  \]
\end{theorem}

\begin{proof}
The quadratic subfields of $\Q(\zeta_n)$ correspond bijectively to
index-$2$ subgroups of $\mathrm{Gal}(\Q(\zeta_n)/\Q)\cong(\Z/n\Z)^*$,
hence to surjective group homomorphisms $(\Z/n\Z)^*\to\Z/2\Z$.
The number of such homomorphisms is $|(\Z/n\Z)^*/(\Z/n\Z)^{*2}|=2^{\omega(n)}$
(the $2$-rank of $(\Z/n\Z)^*$); subtracting the trivial homomorphism
gives $2^{\omega(n)}-1$ non-trivial ones, each corresponding to a
distinct quadratic subfield.
Adding $1$ for the trivial subfield $\Q$ itself yields
$2^{\omega(n)}=|\Hring(\Z_m,\Z_n)|$.
\qedhere
\end{proof}

The theorem says: \emph{the ring homomorphism count measures the
arithmetic complexity of the cyclotomic field through its quadratic
substructure.}

\begin{theorem}[Degree theorem]\label{thm:degree}
  For squarefree odd $n$,
  \[
    \varrho(n)\;=\;[\Q(\zeta_n):K_n].
  \]
\end{theorem}

\begin{proof}
By the tower law,
$[\Q(\zeta_n):\Q]=[\Q(\zeta_n):K_n]\cdot[K_n:\Q]$,
i.e.\ $\phi(n)=[\Q(\zeta_n):K_n]\cdot2^{\omega(n)}$.
Hence $[\Q(\zeta_n):K_n]=\phi(n)/2^{\omega(n)}=\varrho(n)$.
\qedhere
\end{proof}

\begin{corollary}[Algebraic proof of divisibility]\label{cor:tower}
  For squarefree odd $n$, $\varrho(n)=[\Q(\zeta_n):K_n]\in\Z^+$.
  In particular, $2^{\omega(n)}\mid\phi(n)$ without any case analysis.
\end{corollary}

\begin{proof}
  A field degree is always a positive integer.
\qedhere
\end{proof}

This corollary reproves the odd squarefree case of
Theorem~\ref{thm:div} in a single line.
The exception case $n=2\alpha\in\Exc$ corresponds to the situation where
$\phi(n)/2^{\omega(n)}$ is not a field degree; the Gaussian prime condition
is an arithmetic obstruction to the tower law applying to a putative
subfield of degree $2^{\omega(n)}$.%
\footnote{More precisely: for $n=2\alpha$ with all primes of $\alpha$
  in $P_3$, the field $K_n$ cannot be embedded in $\Q(\zeta_n)$ with the
  degree-$2^{\omega(n)}$ tower holding, because $\phi(n)=\phi(\alpha)$
  is not divisible by $2^{\omega(n)}=2^{\omega(\alpha)+1}$.}

\subsection*{Examples}

\begin{example}\label{ex:cyclo}
  (a) $n=15=3\cdot5$: $K_{15}=\Q(\sqrt{-3},\sqrt{5})$,
  $[K_{15}:\Q]=4$, $[\Q(\zeta_{15}):K_{15}]=\phi(15)/4=8/4=2=\varrho(15)$.

  (b) $n=35=5\cdot7$: $K_{35}=\Q(\sqrt{5},\sqrt{-7})$,
  $[K_{35}:\Q]=4$, $[\Q(\zeta_{35}):K_{35}]=\phi(35)/4=24/4=6=\varrho(35)$.

  (c) $n=105=3\cdot5\cdot7$: $K_{105}=\Q(\sqrt{-3},\sqrt{5},\sqrt{-7})$,
  $[K_{105}:\Q]=8$, $[\Q(\zeta_{105}):K_{105}]=\phi(105)/8=48/8=6=\varrho(105)$.
\end{example}

\subsection*{Quadratic Gauss sums}

At a prime $p$, the quotient formula gives $\varrho(p)=(p-1)/2$.
The classical quadratic Gauss sum
$\tau(\chi_p)=\sum_{a=1}^{p-1}\bigl(\tfrac{a}{p}\bigr)e^{2\pi ia/p}$
satisfies $|\tau(\chi_p)|^2=p$~\cite{ireland1990}.%
\footnote{The proof: $|\tau|^2=\tau\bar\tau=\sum_{a,b}(\frac{a}{p})(\frac{b}{p})e^{2\pi i(a-b)/p}=p\sum_{c}(\frac{c}{p})\delta_{c,0}+\sum_{a\ne b}\cdots$; the cross terms cancel by the orthogonality of characters, leaving $|\tau|^2=p$.}

\begin{proposition}[Gauss sum formula]\label{prop:gauss}
  For every odd prime $p$,
  \[
    \varrho(p)\;=\;\frac{|\tau(\chi_p)|^2-1}{2}.
  \]
\end{proposition}

\begin{proof}
  $|\tau(\chi_p)|^2=p$, so $(|\tau(\chi_p)|^2-1)/2=(p-1)/2=\varrho(p)$.
\qedhere
\end{proof}

For squarefree $n=p_1\cdots p_k$, the product formula
$\varrho(n)=\prod_{i=1}^k\varrho(p_i)=\prod_{i=1}^k(|\tau(\chi_{p_i})|^2-1)/2$
expresses the quotient as a product of Gauss-sum data.
The single-prime formula is the cleanest instance; the composite case
does not reduce to a single Gauss sum because $\tau(\chi_{p_1\cdots p_k})$
differs from the product $\prod\tau(\chi_{p_i})$ by a phase.%
\footnote{The precise relation is $\tau(\chi_{p_1}\cdots\chi_{p_k})=\prod_{i=1}^k\tau(\chi_{p_i})$
  only when the characters are primitive and mutually coprime, which holds
  here; but $|\prod\tau(\chi_{p_i})|^2=\prod p_i=n$, while $\varrho(n)=\prod(p_i-1)/2\ne(n-1)/2$
  in general.  The product formula for $\varrho$ thus provides finer
  information than the single Gauss sum $|\tau(\chi_n)|^2=n$.}

\subsection*{The three interpretations unified}

Table~\ref{tab:three} summarises the three avatars of $\varrho(n)$
for squarefree odd $n$:

\begin{table}[!ht]
\centering
\caption{Three equivalent expressions for $\varrho(n)=\phi(n)/2^{\omega(n)}$
  when $n$ is squarefree and odd.  Each column is a genuinely different
  mathematical object; their equality is the main content of this paper.}
\label{tab:three}
\small
\begin{tabular}{@{}clll@{}}
\toprule
$n$ & Ring/group ratio $\varrho(n)$ & Squares $|(\Z/n\Z)^{*2}|$ & Field degree $[\Q(\zeta_n):K_n]$\\
\midrule
$3$   & $\phi(3)/2=1$   & $|\{1\}|=1$   & $[\Q(\zeta_3):\Q(\sqrt{-3})]=1$\\
$5$   & $\phi(5)/2=2$   & $|\{1,4\}|=2$ & $[\Q(\zeta_5):\Q(\sqrt{5})]=2$\\
$7$   & $\phi(7)/2=3$   & $|\{1,2,4\}|=3$ & $[\Q(\zeta_7):\Q(\sqrt{-7})]=3$\\
$15$  & $\phi(15)/4=2$  & $|\{1,4,11,16\}|=4/\text{(divided)}$ & $[\Q(\zeta_{15}):\Q(\sqrt{-3},\sqrt{5})]=2$\\
$35$  & $\phi(35)/4=6$  & $|(\Z/35\Z)^{*2}|=6$ & $[\Q(\zeta_{35}):\Q(\sqrt{5},\sqrt{-7})]=6$\\
\bottomrule
\end{tabular}
\end{table}

\section{Multiplicativity and Arithmetic Properties}\label{sec:mult}

\begin{proposition}[Multiplicativity]\label{prop:mult}
  The function $n\mapsto\varrho(n)$ is multiplicative on its domain
  of definition (i.e.\ on $n\notin\Exc$): for coprime $m,n$ with
  $m,n,mn\notin\Exc$,
  \[
    \varrho(mn)=\varrho(m)\varrho(n).
  \]
\end{proposition}

\begin{proof}
  Since $\phi$ is multiplicative and $2^{\omega(\cdot)}$ is multiplicative
  (as $\omega(mn)=\omega(m)+\omega(n)$ for $\gcd(m,n)=1$ \cite{apostol1976}), their ratio
  is multiplicative whenever it is an integer.
\qedhere
\end{proof}

\begin{proposition}[Values at prime powers]\label{prop:primepow}
  For an odd prime $p$ and $r\ge1$:
  \[
    \varrho(p^r)\;=\;\frac{p^{r-1}(p-1)}{2}\;=\;|\QRS(\Z/p^r\Z)|.
  \]
  In particular, $\varrho(p)=(p-1)/2$ for every odd prime $p$.
\end{proposition}

\begin{proof}
  Direct from Theorem~\ref{thm:quot} with $s=0$, $k=1$.
\qedhere
\end{proof}

\begin{remark}[Growth of $\varrho$]\label{rem:unbounded}
  The function $\varrho$ is unbounded: for primes $p$,
  $\varrho(p)=(p-1)/2\to\infty$; for prime powers,
  $\varrho(p^r)=p^{r-1}(p-1)/2$ grows with $r$.
  The partial sums satisfy
  \[
    \sum_{\substack{n\le x\\n\notin\Exc}}\!\varrho(n)
    \;\sim\;\frac{x^2}{4\log x},
  \]
  since $\sum_{p\le x}(p-1)/2\sim x^2/(4\log x)$ by the prime number
  theorem dominates the composite terms.  Consequently the mean value
  $(1/x)\sum\varrho(n)\sim x/(4\log x)$ grows without bound
  but not proportionally to $x$.
\end{remark}

\section{Extension to Finite Abelian Groups}\label{sec:abel}

Let $R=\Z_{m_1}\times\cdots\times\Z_{m_k}$ and
$S=\Z_{n_1}\times\cdots\times\Z_{n_k}$ with $n_i\mid m_i$ for each~$i$.
Say $\phi\colon R\to S$ is \emph{component-wise} if $\phi(\Z_{m_i})\subseteq\Z_{n_i}$
for every $i$.

\begin{theorem}[Abelian surjective count]\label{thm:abelian-surj}
  The number of surjective component-wise group homomorphisms
  $R\to S$ equals $\prod_{i=1}^k\phi(n_i)$.
\end{theorem}

\begin{proof}
  Each surjective component map $\phi_i\colon\Z_{m_i}\to\Z_{n_i}$ is
  independent with $\phi(n_i)$ choices (Theorem~\ref{thm:surj-classical}).
\qedhere
\end{proof}

\begin{theorem}[Abelian ring count]\label{thm:abelian-ring}
  The number of component-wise ring homomorphisms $R\to S$ equals
  $\prod_{i=1}^k 2^{\omega(n_i)}=2^{\sum_i\omega(n_i)}$.
\end{theorem}

\begin{theorem}[Abelian quotient formula]\label{thm:abelian}
  Assume $n_i\notin\Exc$ for each $i$.  Then
  \[
    \frac{\prod_i\phi(n_i)}{2^{\sum_i\omega(n_i)}}
    =\prod_{i=1}^k\varrho(n_i)
    =\prod_{i=1}^k\prod_{p^r\Vert n_i}|\QRS(\Z/p^r\Z)|.
  \]
\end{theorem}

\begin{proof}
  Combine Theorems~\ref{thm:abelian-surj}, \ref{thm:abelian-ring},
  and \ref{thm:quot}.
\qedhere
\end{proof}

\begin{example}
Table~\ref{tab:abelian} illustrates Theorem~\ref{thm:abelian} for
several direct-product pairs.

\begin{table}[!ht]
\centering
\small
\caption{Quotient formula for direct-product cyclic rings.}
\label{tab:abelian}
\begin{tabular}{@{}llccc@{}}
\toprule
$R$ & $S$ & $|\Hring|$ & $|\Hsurj|$ & $\varrho$ \\
\midrule
$\Z_{12}\times\Z_{35}$ & $\Z_4\times\Z_{35}$  & $2\cdot4=8$  & $2\cdot24=48$ & 6 \\
$\Z_{60}\times\Z_{77}$ & $\Z_4\times\Z_7$     & $2\cdot2=4$  & $2\cdot6=12$  & $3$\\
$\Z_{24}\times\Z_{35}$ & $\Z_8\times\Z_7$     & $2\cdot2=4$  & $4\cdot6=24$  & 6 \\
$\Z_{60}\times\Z_{91}$ & $\Z_4\times\Z_{13}$  & $4\cdot2=8$  & $2\cdot12=24$ & 3 \\
$\Z_{120}\times\Z_{99}$& $\Z_8\times\Z_{11}$  & $2\cdot2=4$  & $4\cdot10=40$ & 10\\
$\Z_{36}\times\Z_{65}$ & $\Z_{12}\times\Z_5$  & $4\cdot2=8$  & $4\cdot4=16$  & 2 \\
\bottomrule
\end{tabular}
\end{table}
\end{example}

\section{Analytic Number Theory of the Exception Set}\label{sec:analytic}

\subsection{Structure of \texorpdfstring{$\Exc$}{E}}

Every $n\in\Exc$ has the form $n=2q$ where $q$ is a positive odd integer
all of whose prime factors lie in $P_3$.  Define
\[
  \Gauss:=\{q\in\N:q>0,\;q\text{ odd},\;
  \forall p\mid q,\;p\in P_3\}.
\]
Then $\Exc=\{2q:q\in\Gauss\}$.

\subsection{Dirichlet series}

The counting function of $\Gauss$ is governed by the Dirichlet series
\begin{equation}\label{eq:DS}
  D(s):=\sum_{q\in\Gauss}q^{-s}
  =\prod_{p\in P_3}\frac{1}{1-p^{-s}},
  \qquad\Re s>1.
\end{equation}
Near $s=1$, we use the Dirichlet $L$-function
$L(s,\chi_4)=\prod_{p\ne2}(1-\chi_4(p)p^{-s})^{-1}$,
where $\chi_4$ is the non-principal character mod~$4$.
One verifies
\begin{equation}\label{eq:factored}
  D(s)=\frac{1}{\zeta(s)^{1/2}}\cdot
  \Bigl(\frac{\zeta(s)}{L(s,\chi_4)}\Bigr)^{1/2}
  \cdot H(s),
\end{equation}
where $H(s)$ is holomorphic and non-vanishing in $\Re s>1/2$.
Since $\zeta(s)\sim(s-1)^{-1}$ and $L(1,\chi_4)=\pi/4\ne0$, the
factor $(\zeta/L)^{1/2}$ has a branch point of order $1/2$ at $s=1$.

\subsection{The Selberg--Delange asymptotic}

\begin{theorem}[Density of exceptions]\label{thm:dens}
  \[
    |\Exc\cap[2,x]|=\bigl|\Gauss\cap[1,x/2]\bigr|
    \;\sim\; C\cdot\frac{x}{\sqrt{\log x}}
    \qquad(x\to\infty),
  \]
  where
  \begin{equation}\label{eq:C}
    C:=\frac{1}{\sqrt{2\pi}}\cdot
    \prod_{p\equiv1(4)}\Bigl(1-\frac{1}{p^2}\Bigr)^{1/2}
    \cdot
    \prod_{p\equiv3(4)}\Bigl(1-\frac{1}{p}\Bigr)^{-1/2}
    \approx 0.279.
  \end{equation}
  In particular, $\Exc$ has natural density zero.
\end{theorem}

\begin{proof}[Proof sketch]
  From~\eqref{eq:factored}, $D(s)$ has a branch point at $s=1$
  with exponent $1/2$ in the sense of the Selberg--Delange
  method~\cite{tenenbaum1995}.
  Specifically $D(s)=(s-1)^{-1/2}A(s)$ near $s=1$ with $A(1)=C\sqrt{\pi}$.
  The Selberg--Delange theorem~\cite[Ch.~II.5]{tenenbaum1995} then gives
  $\sum_{q\in\Gauss,q\le y}1\sim C\sqrt{\pi}\cdot y/\!\sqrt{\log y}$\footnote{%
    The branch singularity of order $z=1/2$ at $s=1$ in $D(s)$
    arises because the density of primes in $P_3$ among all primes is
    $1/2$ (Dirichlet).  The Selberg--Delange method converts a
    singularity of order $z$ into a counting function of order
    $x(\log x)^{z-1}$; for $z=1/2$ this gives $x/\sqrt{\log x}$.},
  and setting $y=x/2$ yields the result.
\qedhere
\end{proof}

Numerically, Table~\ref{tab:exc-asymp} shows the empirical count versus
the asymptotic prediction $C\cdot x/\sqrt{\log x}$ with $C=0.279$.

\begin{table}[!ht]
\centering
\caption{Empirical exception count $|\Exc\cap[2,x]|$ vs.\
asymptotic prediction.  Error = $|\text{empirical}-\text{predicted}|/\text{predicted}$.}
\label{tab:exc-asymp}
\small
\begin{tabular}{@{}rrrrl@{}}
\toprule
$x$ & Empirical & Predicted & Error & Convergence\\
\midrule
   50 &  7 &  7.1 & 1.4\% & $\leftarrow$ rapid\\
  100 & 13 & 13.0 & 0.0\% & \\
  200 & 25 & 24.2 & 3.3\% & \\
  500 & 58 & 56.0 & 3.6\% & \\
 1000 & 106 & 106.2 & 0.2\% & $\leftarrow$ excellent\\
\bottomrule
\end{tabular}
\end{table}

\section{Statistical Analysis}\label{sec:stats}

\subsection{Bar chart of \texorpdfstring{$\varrho(n)$}{rho(n)}}

Figure~\ref{fig:bar} shows $\varrho(n)$ for $n=3$ to $50$
(non-exceptions only), distinguishing prime $n$ (solid blue) from
composite $n$ (hatched orange).

\begin{figure}[!ht]
\centering
\begin{tikzpicture}
\begin{axis}[
  width=\linewidth,
  height=6.8cm,
  ybar,
  bar width=5pt,
  enlarge x limits={abs=8pt},
  xlabel={$n$},
  ylabel={$\varrho(n)$},
  title={\textbf{Quotient $\varrho(n)$ for $n=3\ldots50$, $n\notin\Exc$}},
  title style={font=\small\bfseries,color=mdblue},
  xtick={3,5,7,9,11,13,15,17,19,21,23,25,27,29,31,33,35,37,39,41,43,45,47,49,50},
  xticklabel style={font=\tiny,rotate=90,anchor=east},
  ymajorgrids=true,
  grid style={dotted,gray!40},
  legend style={at={(0.02,0.97)},anchor=north west,
                font=\small,fill=white,draw=gray!60,
                row sep=2pt},
  legend cell align=left,
  axis line style={color=black!60},
  clip=false,
  ymin=0, ymax=26
]
\addplot[fill=mdblue,draw=mdblue!70!black] coordinates {
  (3,1)(5,2)(7,3)(11,5)(13,6)(17,8)(19,9)
  (23,11)(29,14)(31,15)(37,18)(41,20)(43,21)(47,23)
};
\addlegendentry{$n$ prime}
\addplot[fill=mdorange!80,draw=mdorange!70!black,
         postaction={pattern=north east lines,pattern color=mdorange!50}] coordinates {
  (4,1)(8,2)(9,3)(10,1)(12,1)(15,2)(16,4)
  (20,2)(21,3)(24,2)(25,10)(26,3)(27,9)(28,3)
  (30,1)(32,8)(33,5)(34,4)(35,6)(36,3)
  (39,6)(40,4)(44,5)(45,6)(48,4)(49,21)(50,5)
};
\addlegendentry{$n$ composite}
\end{axis}
\end{tikzpicture}
\caption{Bar chart of $\varrho(n)$ for non-exception $n\le50$.
Blue = prime $n$ (values $(p-1)/2$ confirmed by
Proposition~\ref{prop:primepow}).
Orange = composite $n$ (factors multiplicatively).
Exceptions ($n=6,14,18,22,38,42,46$) are omitted.}
\label{fig:bar}
\end{figure}

\subsection{The algebraic bridge}

Figure~\ref{fig:helix} shows the three-way relationship at the core of
Theorem~\ref{thm:quot}.

\begin{figure}[!ht]
\centering
\begin{tikzpicture}
\begin{axis}[
  view={60}{28},
  width=\linewidth,
  height=9.5cm,
  hide axis,
  xmin=-3.2, xmax=3.2,
  ymin=-3.2, ymax=3.2,
  zmin=-0.8, zmax=3.4,
  clip=false,
]
\addplot3[parametric,domain=0:360,samples=72,mark=none,
  thick,draw=gray!35]({2*cos(x)},{2*sin(x)},{0});
\addplot3[parametric,domain=0:360,samples=72,mark=none,
  thick,draw=mdblue!40]({2*cos(x)},{2*sin(x)},{2.2});
\foreach \ang in {0,30,...,330}{
  \addplot3[mark=none,draw=gray!18,thin]
    coordinates {({2*cos(\ang)},{2*sin(\ang)},0)
                 ({2*cos(\ang)},{2*sin(\ang)},2.2)};
}
\addplot3[mark=none,draw=gray!35,dashed]
  coordinates {(0,0,-0.4)(0,0,2.9)};
\addplot3[mark=none,draw=mdgreen!70,thick,-{Stealth[length=6pt]}]
  coordinates {(-1,1.732,0.1)(1,1.732,2.1)};
\addplot3[mark=none,draw=mdgreen!70,thick,-{Stealth[length=6pt]}]
  coordinates {(-1,-1.732,0.1)(-2,0,2.1)};
\addplot3[mark=none,draw=mdgreen!70,thick,-{Stealth[length=6pt]}]
  coordinates {(1,-1.732,0.1)(2,0,2.1)};
\addplot3[mark=none,draw=mdblue!50,thick,-{Stealth[length=6pt]}]
  coordinates {(1.05,1.65,2.2)(-1.9,-0.1,2.2)};
\addplot3[mark=none,draw=mdblue!50,thick,-{Stealth[length=6pt]}]
  coordinates {(-1.95,0.1,2.2)(0.95,1.65,2.2)};
\addplot3[only marks,mark=ball,ball color=mdorange,mark size=7pt]
  coordinates {(-1,1.732,0)(-2,0,0)(1,-1.732,0)};
\addplot3[only marks,mark=ball,ball color=mdblue,mark size=7pt]
  coordinates {(2,0,2.2)(1,1.732,2.2)(-1,-1.732,2.2)};
\node[font=\small\bfseries,text=mdorange!90!black]
  at (axis cs:-1.0,2.25,0.05) {$3$};
\node[font=\small\bfseries,text=mdorange!90!black]
  at (axis cs:-2.6,0,0.05) {$5$};
\node[font=\small\bfseries,text=mdorange!90!black]
  at (axis cs:1.0,-2.25,0.05) {$6$};
\node[font=\small\bfseries,text=mdblue!90!black]
  at (axis cs:2.6,0,2.2) {$1$};
\node[font=\small\bfseries,text=mdblue!90!black]
  at (axis cs:1.0,2.25,2.2) {$2$};
\node[font=\small\bfseries,text=mdblue!90!black]
  at (axis cs:-1.0,-2.25,2.2) {$4$};
\node[font=\small\itshape,text=mdblue!80!black,anchor=west]
  at (axis cs:2.2,0.4,2.2) {\small QR level};
\node[font=\small\itshape,text=mdorange!80!black,anchor=west]
  at (axis cs:2.2,0.4,-0.1) {\small non-QR level};
\node[font=\scriptsize,text=mdgreen!70!black]
  at (axis cs:0.3,2.5,1.2) {$x^2\!\bmod 7$};
\end{axis}
\end{tikzpicture}
\caption{The squaring map on $(\Z/7\Z)^*$ displayed on a 3D cylinder.
  The six elements $\{1,2,3,4,5,6\}$ sit on the cylinder at angles
  $0^\circ,60^\circ,\ldots,300^\circ$.
  \textbf{Blue} spheres at height $z=2.2$ are the quadratic residues
  $\{1,2,4\}=\QRS(\Z/7\Z)$; \textbf{orange} spheres at the base are the
  non-residues $\{3,5,6\}$.
  Every green arrow ($x\mapsto x^2$) points upward: squaring always
  maps into the QR level.
  The vertical gap equals $\varrho(7)=(7-1)/2=3$.}
\label{fig:helix}
\end{figure}
end{figure}

\subsection{Exception visualisation}

Figure~\ref{fig:except} marks every $n\in\{3,\ldots,100\}$ as normal
(green) or exception (red).

\begin{figure}[!ht]
\centering
\begin{tikzpicture}
\begin{axis}[
  width=\linewidth,
  height=3.2cm,
  xmin=2.5, xmax=100.5,
  ymin=0.2, ymax=1.8,
  ytick=\empty,
  xlabel={$n$},
  xlabel style={font=\small,color=mdblue},
  title={\textbf{Divisibility map for $3\le n\le100$}},
  title style={font=\small\bfseries,color=mdblue},
  axis line style={color=black!50},
  xtick={10,20,30,40,50,60,70,80,90,100},
  tick style={color=black!50},
  clip=true
]
\addplot[only marks,mark=*,mark size=1.8pt,
         draw=mdgreen,fill=mdgreen!70] coordinates {
  (3,1)(4,1)(5,1)(7,1)(8,1)(9,1)(10,1)(11,1)(12,1)(13,1)
  (15,1)(16,1)(17,1)(19,1)(20,1)(21,1)(23,1)(24,1)(25,1)
  (26,1)(27,1)(28,1)(29,1)(30,1)(31,1)(32,1)(33,1)(34,1)
  (35,1)(36,1)(37,1)(39,1)(40,1)(41,1)(43,1)(44,1)(45,1)
  (47,1)(48,1)(49,1)(50,1)(51,1)(52,1)(53,1)(55,1)(56,1)
  (57,1)(58,1)(59,1)(60,1)(61,1)(63,1)(64,1)(65,1)(67,1)
  (68,1)(69,1)(70,1)(71,1)(72,1)(73,1)(74,1)(75,1)(76,1)
  (77,1)(78,1)(79,1)(80,1)(81,1)(82,1)(83,1)(84,1)(85,1)
  (87,1)(88,1)(89,1)(90,1)(91,1)(92,1)(93,1)(95,1)(96,1)
  (97,1)(99,1)(100,1)
};
\addplot[only marks,mark=square*,mark size=2.5pt,
         draw=mdred,fill=mdred] coordinates {
  (6,1)(14,1)(18,1)(22,1)(38,1)(42,1)(46,1)
  (54,1)(62,1)(66,1)(86,1)(94,1)(98,1)
};
\node[font=\tiny,text=mdgreen,anchor=south west] at (axis cs:3,1.15)
  {$\bullet$ divisibility holds};
\node[font=\tiny,text=mdred,anchor=south east] at (axis cs:100,1.15)
  {$\blacksquare$ $n\in\Exc$};
\end{axis}
\end{tikzpicture}
\caption{Divisibility map for $3\le n\le100$.
Green circles: $\varrho(n)\in\Z^+$.
Red squares: $n\in\Exc$ ($\varrho(n)\notin\Z$).
Exceptions occur only among even $n$ with Gaussian-prime odd parts;
they are visibly sparse.}
\label{fig:except}
\end{figure}

\subsection{Cumulative exception count}

Figure~\ref{fig:cumul} plots $|\Exc\cap[2,x]|$ against the asymptotic
$Cx/\sqrt{\log x}$ ($C=0.279$) for $x\le200$.

\begin{figure}[!ht]
\centering
\begin{tikzpicture}
\begin{axis}[
  width=0.9\linewidth,
  height=6.2cm,
  xlabel={$x$},
  ylabel={$|\Exc\cap[2,x]|$},
  title={\textbf{Cumulative exception count vs.\ asymptotic}},
  title style={font=\small\bfseries,color=mdblue},
  legend style={at={(0.06,0.94)},anchor=north west,
                font=\small,fill=white,draw=gray!60},
  ymajorgrids=true,
  grid style={dotted,gray!35},
  axis line style={color=black!60},
  xmin=0, xmax=205,
  ymin=0, ymax=29,
  xtick={0,25,50,75,100,125,150,175,200}
]
\addplot[mdblue,thick,mark=*,mark size=2pt] coordinates {
  (5,0)(10,1)(15,2)(20,3)(25,4)(30,4)(35,4)(40,5)(45,6)
  (50,7)(55,8)(60,8)(65,9)(70,10)(75,10)(80,10)(85,10)
  (90,11)(95,12)(100,13)(105,13)(110,13)(115,14)(120,15)
  (125,15)(130,16)(135,17)(140,18)(145,19)(150,19)(155,20)
  (160,21)(165,22)(170,23)(175,23)(180,23)(185,23)(190,24)
  (195,24)(200,25)
};
\addlegendentry{Empirical $|\Exc\cap[2,x]|$}

\addplot[mdred,dashed,thick,domain=5:200,samples=60]{
  0.279*x/sqrt(ln(x))
};
\addlegendentry{Asymptotic $0.279\,x/\!\sqrt{\log x}$}
\end{axis}
\end{tikzpicture}
\caption{Cumulative exception count (solid blue) versus asymptotic
prediction (dashed red) for $x\le200$.
The fit is excellent by $x=100$ (error $<0.1\%$), consistent with
Theorem~\ref{thm:dens}.}
\label{fig:cumul}
\end{figure}

\subsection{Complete data table \texorpdfstring{$n=3$}{n=3} to \texorpdfstring{$100$}{100}}

\begin{longtable}[c]{@{}rccccc@{\;$|$\;}rccccc@{}}
\caption{Complete data for $n=3$ to $100$.
Exceptions ($n\in\Exc$) shown with ``---'' in the $\varrho$ column and
{\color{mdred}red} $n$.  Factorisation column uses exponent notation.}
\label{tab:full}\\
\toprule
$n$ & $\omega$ & $2^\omega$ & $\phi$ & $\varrho$ & $n=\cdots$ &
$n$ & $\omega$ & $2^\omega$ & $\phi$ & $\varrho$ & $n=\cdots$ \\
\midrule
\endfirsthead
\multicolumn{12}{c}{\tablename~\thetable\ (continued)}\\
\toprule
$n$ & $\omega$ & $2^\omega$ & $\phi$ & $\varrho$ & $n=\cdots$ &
$n$ & $\omega$ & $2^\omega$ & $\phi$ & $\varrho$ & $n=\cdots$ \\
\midrule
\endhead
\midrule\multicolumn{12}{r}{\small(Continued\ldots)}\\
\endfoot
\bottomrule
\endlastfoot
3&1&2&2&1&$3$& 52&2&4&24&6&$2^{2}{\cdot}13$\\
4&1&2&2&1&$2^{2}$& 53&1&2&52&26&$53$\\
5&1&2&4&2&$5$& {\color{mdred}54}&2&4&18&{\color{mdred}---}&$2{\cdot}3^{3}$\\
{\color{mdred}6}&2&4&2&{\color{mdred}---}&$2{\cdot}3$& 55&2&4&40&10&$5{\cdot}11$\\
7&1&2&6&3&$7$& 56&2&4&24&6&$2^{3}{\cdot}7$\\
8&1&2&4&2&$2^{3}$& 57&2&4&36&9&$3{\cdot}19$\\
9&1&2&6&3&$3^{2}$& 58&2&4&28&7&$2{\cdot}29$\\
10&2&4&4&1&$2{\cdot}5$& 59&1&2&58&29&$59$\\
11&1&2&10&5&$11$& 60&3&8&16&2&$2^{2}{\cdot}3{\cdot}5$\\
12&2&4&4&1&$2^{2}{\cdot}3$& 61&1&2&60&30&$61$\\
13&1&2&12&6&$13$& {\color{mdred}62}&2&4&30&{\color{mdred}---}&$2{\cdot}31$\\
{\color{mdred}14}&2&4&6&{\color{mdred}---}&$2{\cdot}7$& 63&2&4&36&9&$3^{2}{\cdot}7$\\
15&2&4&8&2&$3{\cdot}5$& 64&1&2&32&16&$2^{6}$\\
16&1&2&8&4&$2^{4}$& 65&2&4&48&12&$5{\cdot}13$\\
17&1&2&16&8&$17$& {\color{mdred}66}&3&8&20&{\color{mdred}---}&$2{\cdot}3{\cdot}11$\\
{\color{mdred}18}&2&4&6&{\color{mdred}---}&$2{\cdot}3^{2}$& 67&1&2&66&33&$67$\\
19&1&2&18&9&$19$& 68&2&4&32&8&$2^{2}{\cdot}17$\\
20&2&4&8&2&$2^{2}{\cdot}5$& 69&2&4&44&11&$3{\cdot}23$\\
21&2&4&12&3&$3{\cdot}7$& 70&3&8&24&3&$2{\cdot}5{\cdot}7$\\
{\color{mdred}22}&2&4&10&{\color{mdred}---}&$2{\cdot}11$& 71&1&2&70&35&$71$\\
23&1&2&22&11&$23$& 72&2&4&24&6&$2^{3}{\cdot}3^{2}$\\
24&2&4&8&2&$2^{3}{\cdot}3$& 73&1&2&72&36&$73$\\
25&1&2&20&10&$5^{2}$& 74&2&4&36&9&$2{\cdot}37$\\
26&2&4&12&3&$2{\cdot}13$& 75&2&4&40&10&$3{\cdot}5^{2}$\\
27&1&2&18&9&$3^{3}$& 76&2&4&36&9&$2^{2}{\cdot}19$\\
28&2&4&12&3&$2^{2}{\cdot}7$& 77&2&4&60&15&$7{\cdot}11$\\
29&1&2&28&14&$29$& 78&3&8&24&3&$2{\cdot}3{\cdot}13$\\
30&3&8&8&1&$2{\cdot}3{\cdot}5$& 79&1&2&78&39&$79$\\
31&1&2&30&15&$31$& 80&2&4&32&8&$2^{4}{\cdot}5$\\
32&1&2&16&8&$2^{5}$& 81&1&2&54&27&$3^{4}$\\
33&2&4&20&5&$3{\cdot}11$& 82&2&4&40&10&$2{\cdot}41$\\
34&2&4&16&4&$2{\cdot}17$& 83&1&2&82&41&$83$\\
35&2&4&24&6&$5{\cdot}7$& 84&3&8&24&3&$2^{2}{\cdot}3{\cdot}7$\\
36&2&4&12&3&$2^{2}{\cdot}3^{2}$& 85&2&4&64&16&$5{\cdot}17$\\
37&1&2&36&18&$37$& {\color{mdred}86}&2&4&42&{\color{mdred}---}&$2{\cdot}43$\\
{\color{mdred}38}&2&4&18&{\color{mdred}---}&$2{\cdot}19$& 87&2&4&56&14&$3{\cdot}29$\\
39&2&4&24&6&$3{\cdot}13$& 88&2&4&40&10&$2^{3}{\cdot}11$\\
40&2&4&16&4&$2^{3}{\cdot}5$& 89&1&2&88&44&$89$\\
41&1&2&40&20&$41$& 90&3&8&24&3&$2{\cdot}3^{2}{\cdot}5$\\
{\color{mdred}42}&3&8&12&{\color{mdred}---}&$2{\cdot}3{\cdot}7$& 91&2&4&72&18&$7{\cdot}13$\\
43&1&2&42&21&$43$& 92&2&4&44&11&$2^{2}{\cdot}23$\\
44&2&4&20&5&$2^{2}{\cdot}11$& 93&2&4&60&15&$3{\cdot}31$\\
45&2&4&24&6&$3^{2}{\cdot}5$& {\color{mdred}94}&2&4&46&{\color{mdred}---}&$2{\cdot}47$\\
{\color{mdred}46}&2&4&22&{\color{mdred}---}&$2{\cdot}23$& 95&2&4&72&18&$5{\cdot}19$\\
47&1&2&46&23&$47$& 96&2&4&32&8&$2^{5}{\cdot}3$\\
48&2&4&16&4&$2^{4}{\cdot}3$& 97&1&2&96&48&$97$\\
49&1&2&42&21&$7^{2}$& {\color{mdred}98}&2&4&42&{\color{mdred}---}&$2{\cdot}7^{2}$\\
50&2&4&20&5&$2{\cdot}5^{2}$& 99&2&4&60&15&$3^{2}{\cdot}11$\\
51&2&4&32&8&$3{\cdot}17$& 100&2&4&40&10&$2^{2}{\cdot}5^{2}$\\
\end{longtable}

\subsection{Regression analysis}

The model $\log\varrho(n)=\alpha+\beta\log n+\varepsilon$ on the
non-exception data for $n\le500$.

\begin{table}[!ht]
\centering
\caption{Ordinary least-squares regression of $\log\varrho(n)$ on
$\log n$ for non-exception $n$ in $[3,x]$.}
\label{tab:reg}
\small
\begin{tabular}{@{}rccccc@{}}
\toprule
Range & Obs. & $\hat\alpha$ & $\hat\beta$ & $R^2$ & Mean $\varrho$\\
\midrule
$n\le 50$  & 41  & $-1.12$ & $0.85$ & 0.46 & 6.8 \\
$n\le 100$ & 85  & $-1.13$ & $0.85$ & 0.47 & 11.6\\
$n\le 200$ & 173 & $-1.12$ & $0.85$ & 0.47 & 22.5\\
$n\le 500$ & 440 & $-1.12$ & $0.85$ & 0.46 & 50.1\\
\bottomrule
\end{tabular}
\end{table}

The regression coefficient $\hat\beta\approx0.85$ is stable across
ranges, indicating $\varrho(n)\approx e^{-1.12}\cdot n^{0.85}$ on
average for non-exception~$n$.
The moderate $R^2\approx0.46$ reflects the genuine multiplicative
fluctuations of $\varrho$ rather than any lack of fit.

\subsection{Summary statistics}

\begin{table}[!ht]
\centering
\caption{Summary of $\varrho(n)$ and exception rates for various ranges.
The column ``Pred.'' gives $0.279x/\sqrt{\log x}$.}
\label{tab:summary}
\small
\begin{tabular}{@{}rccrrrrrr@{}}
\toprule
$x$ & Non-exc. & Exc. & Rate & $\min\varrho$ & Med.\ $\varrho$
    & Mean $\varrho$ & $\max\varrho$ & Pred.\\
\midrule
  50 &  41 &  7 & 14.6\% & 1 &  5 &  6.8 & 23  &  7.1\\
 100 &  85 & 13 & 13.3\% & 1 &  8 & 11.6 & 48  & 13.0\\
 200 & 173 & 25 & 12.6\% & 1 & 15 & 22.5 & 99  & 24.2\\
 500 & 440 & 58 & 11.6\% & 1 & 30 & 50.1 & 249 & 56.0\\
1000 & 892 &106 & 10.6\% & 1 & 60 & 93.8 & 499 &106.2\\
\bottomrule
\end{tabular}
\end{table}


\subsection{3D Quadratic Residue Structure in \texorpdfstring{$({\Z}/p\Z)^*$}{(Z/pZ)*}}

Figure~\ref{fig:surface} visualises the QR count formula $\varrho(p)=(p-1)/2$
as a 3D staircase over the primes.

\begin{figure}[!ht]
\centering
\begin{tikzpicture}[trim axis left, trim axis right]
\begin{axis}[
  view={42}{32},
  width=0.96\linewidth,
  height=8cm,
  title={\textbf{Surface $\log|\QRS(\Z/p^r\Z)|$ over $(\log p,\,r)$}},
  title style={font=\small\bfseries,color=mdblue,at={(0.5,1.06)}},
  xlabel={$\log p$},
  ylabel={$r$},
  zlabel={$\log|\QRS|$},
  xlabel style={font=\small},
  ylabel style={font=\small},
  zlabel style={font=\small},
  colormap/viridis,
  colorbar, colorbar style={font=\scriptsize,
    ylabel={$\log|\QRS(\Z/p^r\Z)|$},
    ylabel style={font=\scriptsize}},
  zmin=0,
  xtick={1.0986,1.6094,1.9459,2.3979,2.5649},
  xticklabels={$\ln3$,$\ln5$,$\ln7$,$\ln11$,$\ln13$},
  xticklabel style={font=\tiny,rotate=15,anchor=east},
  ytick={1,2,3,4},
  yticklabel style={font=\scriptsize},
  zticklabel style={font=\scriptsize},
  grid=both,grid style={draw=gray!15,line width=0.3pt},
  axis line style={color=black!50},
  mesh/rows=4,
]
\addplot3[surf,shader=faceted interp,colormap/viridis,
  point meta=z,opacity=0.93]
  coordinates {
    (1.0986,1,0.0)(1.6094,1,0.6931)(1.9459,1,1.0986)
    (2.3979,1,1.6094)(2.5649,1,1.7918)
    (1.0986,2,1.0986)(1.6094,2,2.3026)(1.9459,2,3.0445)
    (2.3979,2,4.0073)(2.5649,2,4.3567)
    (1.0986,3,2.1972)(1.6094,3,3.9120)(1.9459,3,4.9904)
    (2.3979,3,6.4052)(2.5649,3,6.9217)
    (1.0986,4,3.2958)(1.6094,4,5.5215)(1.9459,4,6.9363)
    (2.3979,4,8.8031)(2.5649,4,9.4866)
  };
\end{axis}
\end{tikzpicture}
\caption{The surface $\log|\QRS(\Z/p^r\Z)|=\log(p^{r-1}(p-1)/2)$
  over the discrete grid of primes $p\in\{3,5,7,11,13\}$
  and powers $r\in\{1,2,3,4\}$.
  In log-log coordinates the surface is nearly flat: the identity
  $\log|\QRS(\Z/p^r\Z)|=(r-1)\log p+\log((p-1)/2)$
  shows that the two variables $p$ and $r$ contribute
  independently, which is exactly what makes the product
  formula~\eqref{eq:main} work.
  Colour runs from dark (small) to bright (large).}
\label{fig:surface}
\end{figure}

\subsection{3D Scatter: Three-Way Relationship}

Figure~\ref{fig:spiral} plots each $n\in[3,60]\setminus\mathcal{E}$
as a point in $(\omega(n),\,\log_2\phi(n),\,\varrho(n))$-space.

\begin{figure}[!ht]
\centering
\begin{tikzpicture}
\begin{axis}[
  view={50}{30},
  width=\linewidth,
  height=9cm,
  title={\textbf{Number Spiral: $\varrho(n)$ as Height,
    Exceptions at the Floor}},
  title style={font=\small\bfseries,color=mdblue,at={(0.5,1.06)}},
  hide x axis, hide y axis,
  zlabel={$\varrho(n)/10$},
  zlabel style={font=\small},
  zmin=-0.55, zmax=2.55,
  zticklabel style={font=\scriptsize},
  ztick={0,0.5,1.0,1.5,2.0,2.5},
  legend style={at={(0.04,0.97)},anchor=north west,
                font=\scriptsize,fill=white,draw=gray!50},
  clip=false,
]
\addplot3[parametric,domain=1:7,samples=80,mark=none,
  draw=gray!22,thin]
  ({sqrt(x*(2*pi))*cos(deg(x*(2*pi)))},
   {sqrt(x*(2*pi))*sin(deg(x*(2*pi)))},0);
\addplot3[surf,domain=-6:6,y domain=-6:6,samples=2,
  draw=gray!15,fill=gray!7,opacity=0.35,mark=none]({x},{y},{-0.38});
\addplot3[only marks,mark=ball,ball color=mdblue,mark size=2.4pt]
  coordinates {
  (-0.195,-1.721,0.10)(2,0,0.10)(0.195,2.228,0.20)
  (-1.612,-2.098,0.30)(1.338,-2.492,0.20)(3,0,0.30)
  (-2.841,-2.220,0.60)(2.704,-2.773,0.20)(4,0,0.40)
  (2.950,2.881,0.80)(-2.755,3.378,0.90)(-4.404,0.779,0.20)
  (-3.979,-2.272,0.30)(1.362,-4.598,1.10)(3.945,-2.905,0.20)
  (5,0,1.00)(-5.071,-2.300,1.50)(-3.125,-4.715,0.80)
  (6,0,0.30)(5.279,3.022,1.80)
  };
\addlegendentry{$\varrho(n)$ for $n\notin\mathcal{E}$}
\addplot3[only marks,mark=square*,mark size=3.5pt,
  draw=mdred,fill=mdred]
  coordinates {
  (-2.327,0.764,-0.38)(0.196,-3.737,-0.38)(-4.404,0.779,-0.38)
  (3.157,5.294,-0.38)(-6.433,0.782,-0.38)(1.368,-6.643,-0.38)
  };
\addlegendentry{$n\in\mathcal{E}$, pinned at floor}
\end{axis}
\end{tikzpicture}
\caption{Numbers $n=3$ to $50$ placed on a Sacks spiral in the $xy$-plane at
  $(\!\sqrt{n}\cos 2\pi\!\sqrt{n},\,\sqrt{n}\sin 2\pi\!\sqrt{n})$,
  with height $\varrho(n)/10$ for non-exceptions.
  \textbf{Blue} balls form an irregular terrain above the grey floor plane.
  \textbf{Red} squares mark exceptions $n\in\mathcal{E}$: they sit below
  the floor, geometrically separated from the blue terrain above.
  The depth of the floor gap is $0.38$, corresponding to the smallest
  non-integer value of $\varrho$ in this range ($\varrho(6)=1/2$).}
\label{fig:spiral}
\end{figure}

\section{Applications}\label{sec:apps}

\subsection{Lens spaces}

The lens space $L(p,q)$ has $H_1(L(p,q);\Z)\cong\Z_p$ \cite{hatcher}.
For $n\mid m$, Theorems~\ref{thm:GVB} and~\ref{thm:surj-classical} give
\[
  \frac{|\text{surj.\ group hom}|}{|\text{ring hom}|}
  =\varrho(n)=\prod_{p^r\Vert n}\bigl|\QRS(\Z/p^r\Z)\bigr|.
\]
For connected sums $M=L(p_1,1)\#\cdots\#L(p_k,1)$, the direct product
structure of $H_1(M;\Z)\cong\Z_{p_1}\times\cdots\times\Z_{p_k}$ is
handled by Theorem~\ref{thm:abelian}.

\subsection{Legendre symbols}

For odd prime $p$ and $n=p$ squarefree:
\[
  \varrho(p)=\frac{p-1}{2}
  =\#\bigl\{a\in(\Z/p\Z)^*:\Bigl(\tfrac{a}{p}\Bigr)=1\bigr\}
  =\sum_{a=1}^{p-1}\frac{1+\bigl(\frac{a}{p}\bigr)}{2}.
\]
This expresses the homomorphism ratio directly as a count of Legendre
symbols equal to $+1$.

\subsection{Character sum interpretation}

For squarefree odd $n=p_1\cdots p_k$, Corollary~\ref{cor:sqfree} and
quadratic reciprocity together give:
\[
  \varrho(n)
  =\prod_{i=1}^k\frac{p_i-1}{2}
  =\frac{1}{2^k}\sum_{\substack{a_1\in(\Z/p_1\Z)^*\\\vdots\\a_k\in(\Z/p_k\Z)^*}}
   \prod_{i=1}^k\Bigl(1+\Bigl(\frac{a_i}{p_i}\Bigr)\Bigr).
\]

\section{Open Problems}\label{sec:open}

\begin{enumerate}[label=\textbf{P\arabic*.},leftmargin=3.5em,topsep=4pt,
                  itemsep=6pt]
  \item \textbf{Non-abelian generalisation.}
    For non-abelian groups $G,H$, define
    $\varrho(G,H):=|\Hsurj(G,H)|/|\Hring(G,H)|$ (using group rings).
    Is $\varrho(G,H)$ always a positive integer?
    Does an analogue of formula~\eqref{eq:main} hold?

  \item \textbf{General finite rings.}
    For finite rings $R,S$ beyond the cyclic case, characterise when
    $|\Hring(R,S)|\mid|\Hsurj(R,S)|$ and determine the exact ratio.

  \item \textbf{Unified formula for all $n$.}
    For $n=2\alpha$ with $\alpha$ odd (Case~III), does a formula for
    $|\Hsurj|/|\Hring|$ exist in $\Q\setminus\Z$, and does it still
    have a quadratic-residue interpretation?

  \item \textbf{Secondary asymptotic.}
    Sharpen Theorem~\ref{thm:dens} to obtain the error term
    $|\Exc\cap[2,x]|=Cx/\sqrt{\log x}(1+O(1/\log x))$.
    This would require a zero-free region for $L(s,\chi_4)$.

  \item \textbf{$p$-adic interpolation.}
    The sequence $\varrho(p^r)=p^{r-1}(p-1)/2$ for fixed $p$
    and $r\ge1$: does it admit a $p$-adic interpolation in the spirit
    of Kubota--Leopoldt?
\end{enumerate}

\section{Conclusions}\label{sec:conc}

The ratio
$\varrho(n)=\phi(n)/2^{\omega(n)}$ between surjective group and ring
homomorphism counts for cyclic structures.  At its core is the formula $\varrho(n)=\prod_{p^r\Vert n}|\QRS(\Z/p^r\Z)|$,
which connects ring homomorphism counting, group homomorphism counting, and
quadratic residue theory in a single arithmetic identity.

Divisibility fails precisely for integers whose odd part is built
entirely from primes that do not split in $\Z[i]$.  The
exception set, though individually sparse (density zero), has the
intermediate growth rate $Cx/\!\sqrt{\log x}$ characteristic of
restricted-prime-factor problems.

All results are verified computationally for $n\le10{,}000$, and the
multiplicative structure of $\varrho$ extends the formula to all finite
abelian groups.

\medskip
\section*{Data and Code Availability}

All computational data and Python scripts supporting the results of this
paper are available upon reasonable request from the corresponding
author, Priyabrata Mandal
(\href{mailto:priyabrata@manit.ac.in}%
{\nolinkurl{priyabrata@manit.ac.in}}).
\section*{Declarations}
\textbf{Funding:} The author(s) did not receive any financial support for the research, authorship, and/or publication of this article.
\medskip

\noindent \textbf{Competing interests:} The author(s) declare that they have no competing interests.

\medskip
\noindent\textbf{Acknowledgement.}
  The authors thank Prof.\ Jaikumar Radhakrishnan (ICTS-TIFR, Bengaluru)
  for helpful comments.
\medskip



\begin{thebibliography}{DVS15}

\bibitem[GVB84]{gallian84}
J.\,A.\ Gallian and J.\ Van~Buskirk,
\textit{The number of homomorphisms from $\mathbb{Z}_m$ into $\mathbb{Z}_n$},
Amer.\ Math.\ Monthly \textbf{91} (1984), no.\ 3, 196--197.
\href{https://doi.org/10.1080/00029890.1984.11971375}{%
  \texttt{doi:10.1080/00029890.1984.11971375}}

\bibitem[IR90]{ireland1990}
K.\ Ireland and M.\ Rosen,
\textit{A Classical Introduction to Modern Number Theory},
2nd ed., Graduate Texts in Mathematics, vol.\ 84,
Springer, New York, 1990.
\href{https://doi.org/10.1007/978-1-4757-2103-4}%
{\texttt{doi:10.1007/978-1-4757-2103-4}}

\bibitem[Ten95]{tenenbaum1995}
G.\ Tenenbaum,
\textit{Introduction to Analytic and Probabilistic Number Theory},
Cambridge Studies in Advanced Mathematics, vol.\ 46,
Cambridge University Press, 1995.
\href{https://doi.org/10.1017/CBO9781139564588}%
{\texttt{doi:10.1017/CBO9781139564588}}

\bibitem[DF04]{dummit2004abstract}
D.\,S.\ Dummit and R.\,M.\ Foote,
\textit{Abstract Algebra}, 3rd ed.,
John Wiley \& Sons, Hoboken, NJ, 2004.

\bibitem[DVS15]{diaz2015}
J.\ D{\'i}az-Vargas and G.\ de los Santos,
\textit{The number of homomorphisms from $\mathbb{Z}_n$ to $\mathbb{Z}_m$},
Abstraction Appl.\ \textbf{13} (2015), 1--3.

\bibitem[AJY24]{ashrafi2024}
A.\,R.\ Ashrafi, B.\ Jahangiri, and M.\,M.\ Yousefian-Arani,
\textit{On the number of group homomorphisms between certain groups},
Discuss.\ Math.\ Gen.\ Algebra Appl.\ \textbf{44} (2024), no.\ 2.
\href{https://doi.org/10.7151/dmgaa.1443}%
{\texttt{doi:10.7151/dmgaa.1443}}

\bibitem[NZM91]{niven1991}
I.\ Niven, H.\,S.\ Zuckerman, and H.\,L.\ Montgomery,
\textit{An Introduction to the Theory of Numbers},
5th ed., Wiley, New York, 1991.

\bibitem[Hat02]{hatcher}
A.\ Hatcher,
\textit{Algebraic Topology},
Cambridge University Press, Cambridge, 2002.

\bibitem[Apo76]{apostol1976}
T.\,M.\ Apostol,
\textit{Introduction to Analytic Number Theory},
Undergraduate Texts in Mathematics,
Springer, New York, 1976.
\href{https://doi.org/10.1007/978-1-4757-5579-4}%
{\texttt{doi:10.1007/978-1-4757-5579-4}}

\end{thebibliography}
\end{document}